\documentclass[12pt]{amsart}
\usepackage{latexsym,amssymb}  
\theoremstyle{plain}
  \newtheorem{thm}{Theorem}[section]
  \newtheorem{prop}[thm]{Proposition}
  \newtheorem{lem}[thm]{Lemma}
  \newtheorem{cor}[thm]{Corollary}
\theoremstyle{definition}
  \newtheorem{dfn}[thm]{Definition}
  \newtheorem{exmp}[thm]{Example}
\theoremstyle{remark}
  \newtheorem{rem}[thm]{Remark}
\numberwithin{equation}{section}
\def\ba{{\bf a}}

\def\bR{{\bf R}}
\def\nlim{\lim\limits_{ \to}}
\def\b1{{\mathbf 1}}

\def\L{T^\bullet}
\def\F{F^\bullet}
\def\M{M^\bullet}
\def\N{N^\bullet}
\def\P{P^\bullet}
\def\Q{Q^\bullet}

\def\I{I^\bullet}

\def\cS{{\mathcal S}}
\def\cE{{\mathcal E}}
\def\DF{{\mathcal F}}
\def\DG{{\mathcal G}}
\def\DFn{{\mathcal F}^*}
\def\DGn{{\mathcal G}^*}
\def\DFo{{\mathcal F}^{\sf op}}
\def\DGo{{\mathcal G}^{\sf op}}
\def\bA{{\mathbf {A}}}
\def\m{{\mathfrak m}}
\def\n{{\mathfrak n}}
\def\NN{{\mathbb N}}
\def\Sq{\operatorname{Sq}(S)}
\def\SqE{\operatorname{Sq}(E)}
\def\Syz{\Omega}
\def\lpd{\operatorname{lpd}}
\def\ZZ{{\mathbb Z}}

\def\image{\operatorname{im}}
\def\reg{\operatorname{reg}}
\def\pdim{\operatorname{proj.dim}}
\def\depth{\operatorname{depth}}
\def\cod{\operatorname{ht}}
\def\CS{{\mathcal S}}

\def\Hom{\operatorname{Hom}}
\def\uHom{\underline{\Hom}}
\def\Ext{\operatorname{Ext}}
\def\uExt{\underline{\Ext}}
\def\plin{\operatorname{proj.lin}}
\def\ilin{\operatorname{inj.lin}}
\def\Tor{\operatorname{Tor}}

\def\coker{\operatorname{coker}}

\def\supp{\operatorname{supp}}
\def\<{{\langle}}
\def\>{{\rangle}}
\def\too{{\longrightarrow}}

\def\Id{\operatorname{Id}}
\def\D{{\mathcal D}^\bullet}
\def\grS{\operatorname{gr} S}

\def\Ao{A^{\sf op}}
\def\Aqo{(A^!)^\op}
\def\grAqo{\operatorname{gr} \Aqo}
\def\GrA{{\operatorname{Gr} A}}
\def\GrAo{\operatorname{Gr} {\Ao}}
\def\grm{\operatorname{gr}_\m}

\def\grA{\operatorname{gr} A}
\def\grAq{\operatorname{gr} {A^!}}
\def\grB{\operatorname{gr} B}
\def\lfA{\operatorname{lf} A}
\def\lfAo{\operatorname{lf} {\Ao}}
\def\grAqo{\operatorname{gr} {(A^!)^{\sf op}}}
\def\op{{\sf op}}
\def\grAo{\operatorname{gr} {\Ao}}

\def\DA{{\bf D}_A}
\def\DAq{{\bf D}_{A^!}}
\def\DAqo{{\bf D}_{{(A^!)^{\sf op}}}}
\def\DS{{\bf D}_S}
\def\DE{{\bf D}_E}
\def\mmEn{\operatorname{*gr} E}
\def\mmSn{\operatorname{*gr} S}
\def\MMSn{\operatorname{*Gr} S}
\def\cH{\mathcal H}
\begin{document}
\title[Castelnuovo-Mumford regularity and weakly Koszul modules]
{Castelnuovo-Mumford regularity for complexes and weakly Koszul modules}
\author{Kohji Yanagawa}
\address{Department of Mathematics, 
Graduate School of Science, Osaka University, Toyonaka, Osaka 
560-0043, Japan}
\email{yanagawa@math.sci.osaka-u.ac.jp}

\begin{abstract}
Let $A$ be a noetherian AS regular Koszul quiver algebra 
(if $A$ is commutative, it is essentially a polynomial ring), 
and $\grA$ the category of finitely generated graded left $A$-modules. 
Following J\o rgensen, we define the Castelnuovo-Mumford regularity 
$\reg(\M)$ of a complex $\M \in D^b(\grA)$ in terms of the local cohomologies 
or the minimal projective resolution of $\M$. 
Let $A^!$ be the  quadratic dual ring of $A$. 
For the Koszul duality functor $\DG : D^b(\grA) \to D^b(\grAq)$, 
we have $\reg(\M) = \max \{ \,  i \mid H^i(\DG (\M)) \ne 0 \, \}$.  
Using these concepts, 
we interpret results of Martinez-Villa and  Zacharia concerning  
{\it weakly Koszul modules} (also called {\it componentwise linear modules}) 
over $A^!$.   As an application, refining a result of 
Herzog and R\"omer, we show that  if $J$ is 
a monomial ideal of an exterior algebra $E=\bigwedge \< y_1, \ldots, y_d \>$, 
$d \geq 3$, then the $(d-2)^{\rm nd}$ syzygy of $E/J$ is weakly Koszul. 
\end{abstract}

\maketitle

\section{Introduction}
Let $S = K[x_1, \ldots, x_d]$ be a polynomial ring over a field $K$. 
We regard $S$ as a graded ring with $\deg x_i = 1$ for all $i$. 
The following is a well-known result. 

\begin{thm}[{c.f. \cite{EG}}]\label{flush point}
Let $M$ be a finitely generated graded $S$-module. 
For an integer $r$, the following conditions are equivalent. 
\begin{itemize}
\item[(1)] $H_\m^i(M)_j = 0$ for all $i, j \in \ZZ$ with $i+j >r$. 
\item[(2)] The truncated module $M_{\geq r} := \bigoplus_{i \geq r} M_i$ 
has an $r$-linear free resolution. 
\end{itemize}
Here $\m := (x_1, \ldots, x_d)$ is the irrelevant ideal of $S$, 
and $H^i_\m(M)$ is the $i^{\rm th}$ local cohomology module. 
\end{thm}

If the conditions of Theorem~\ref{flush point} are satisfied, we say 
$M$ is {\it $r$-regular}. For a sufficiently large $r$, $M$ is $r$-regular.
We call $\reg(M) = \min \{ \,  r \mid  \text{$M$ is $r$-regular} \, \}$ 
the {\it Castelnuovo-Mumford regularity} of $M$. 
This is a very important invariant in commutative algebra.  

Let $A$ be a noetherian AS regular Koszul quiver algebra 
with the graded Jacobson radical $\m := \bigoplus_{i\geq 1} A_i$.  
If $A$ is commutative, $A$ is essentially a polynomial ring. 
When $A$ is connected (i.e.,  $A_0 = K$), 
it is the coordinate ring of a 
``noncommutative projective space" in noncommutative 
algebraic geometry.  Let $\grA$ be the category of finitely generated 
graded left $A$-modules and their degree preserving maps. 
(For a graded ring $B$, $\grB$ means the similar category for $B$.) 
The local cohomology module $H_\m^i(M)$ of 
$M \in \grA$ behaves pretty much like in the commutative case. 
For example, we have ``Serre duality theorem" for the derived category 
$D^b(\grA)$. See \cite{M98, Ye} and 
Theorem~\ref{Serre duality} below. By virtue of this duality, we can show 
that Theorem~\ref{flush point} also holds for {\it bounded complexes} in 
$\grA$. 

\begin{thm}\label{2nd thm}
For a complex $\M \in D^b(\grA)$ and an integer $r$, 
the following conditions are equivalent. 
\begin{itemize}
\item[(1)] $H_\m^i(\M)_j = 0$ for all $i, j \in \ZZ$ with $i+j >r$. 
\item[(2)] The truncated complex $(\M)_{\geq r}$ has an $r$-linear 
projective resolution. 
\end{itemize}
Here $(\M)_{\geq r}$ is the subcomplex of $\M$ whose $i^{\rm th}$ term is 
$(M^i)_{\geq (r-i)}$. 
\end{thm}

For a sufficiently large $r$, the conditions of the above theorem are 
satisfied. The regularity $\reg(\M)$ of $\M$ 
is defined in the natural way. When $A$ is connected, J\o rgensen 
(\cite{Jor04}) has studied the regularity of complexes, and essentially proved 
the above result. See also \cite{Jor99, Mori}. 
(Even in the case when $A$ is a polynomial ring, it seems that 
nobody had considered Theorem~\ref{2nd thm} before \cite{Jor04}.) 
But his motivation and treatment are slightly different from ours. 

For $\M \in D^b(\grA)$, set $\cH(\M)$ to be a complex 
such that $\cH(\M)^i = H^i(M)$ for all $i$ 
and the differential maps are zero. 
Then we have $\reg(\cH(\M)) \geq \reg(\M)$. 
The difference $\reg(\cH(\M)) - \reg(\M)$
is a theme of the last section of this paper. 

Let $A^!$ be the quadratic dual ring of $A$. 
For example, if $S=K[x_1, \ldots, x_d]$ is a polynomial ring, then 
$S^!$ is an exterior algebra $E= \bigwedge \< y_1, \ldots, y_d \>$. 
It is known that $A^!$ is always Koszul, finite dimensional, 
and selfinjective.
The Koszul duality functors 
$\DF : D^b(\grAq) \to D^b(\grA)$ and 
$\DG : D^b(\grA) \to D^b(\grAq)$
give a category equivalence  $D^b(\grAq) \cong D^b(\grA)$ 
(see \cite{BGS}). It is easy to check that 
$$\reg(\M) = \max \{ \, i \mid  H^i(\DG(\M)) \ne 0 \, \}$$
for $\M \in D^b(\grA)$. 

Let $\grAo$ be the category of finitely generated graded {\it right} 
$A$-modules. The above results on $\grA$ also hold for 
$\grAo$. Moreover, we have 
$$\reg( \, \bR \uHom_A(\M, \D) \, ) = 
- \min \{ \, i \mid  H^i(\DG(\M)) \ne 0 \, \}$$
for $\M \in D^b(\grA)$. 
Here $\D$ is a balanced dualizing complex of $A$, which 
gives duality functors between $D^b(\grA)$ and $D^b(\grAo)$. 

Let $B$ be a noetherian Koszul algebra. 
For $M \in \grB$ and $i \in \ZZ$, $M_{\< i \>}$ denotes the 
submodule of $M$ generated by the degree $i$ component $M_i$ of $M$. 
We say $M$ is {\it weakly Koszul} if $M_{\< i \>}$ has a linear projective 
resolution for all $i$. This definition is different 
from the original one given in \cite{MZ}, but they are equivalent. 
(Weakly Koszul modules are also called ``componentwise linear modules" 
by some commutative algebraists.)  
 Martinez-Villa and Zacharia proved that if $N \in \grAq$ then 
the $i^{\rm th}$ syzygy $\Omega_i(N)$
of $N$ is weakly Koszul for $i \gg 0$. For $N \in \grAq$, set 
$$\lpd(N):= \min\{ \, i \in \NN  \mid  \text{$\Syz_i(N)$ is 
weakly Koszul} \, \}.$$  

Let $N \in \grAq$ and  $N' :=\uHom_{A^!}(N, A^!) \in \grAqo$ its dual. 
In Theorem~\ref{wKos}, we show that $N$ is weakly Koszul if and only if 
$\reg(\cH \circ \DFo (N')) = 0$, where 
$\DFo : D^b(\grAqo) \to D^b(\grAo)$ is the Koszul duality functor. 
(Since $\reg(\DFo(N')) = 0$, we have $\reg(\cH \circ \DFo (N')) 
\geq 0$ in general.) Moreover, we have 
$$\lpd(N) = \reg( \, \cH \circ \DFo (N') \,)$$
(Corollary~\ref{wKos=regH}). 
As an application of this formula, we refine a result of 
Herzog and R\"omer on monomial ideals of an exterior algebra. 
Among other things, in Proposition~\ref{d-2}, we show that 
if $J$ is a monomial ideal of an exterior algebra 
$E= \bigwedge \< y_1, \ldots, y_d \>$, $d \geq 3$, then 
$\lpd(E/J) \leq d-2$.

Finally, we remark that Herzog and Iyengar (\cite{HI}) studied the invariant 
$\lpd$ and related concepts over noetherian  commutative (graded) 
local rings. Among other things, they proved that 
$\lpd(N)$ is always finite over some ``nice" local 
rings (e.g., complete intersections whose associated graded rings 
are Koszul).   

\section{Preliminaries}
Let $K$ be a field. 
The ring $A$ treated in this paper is a (not necessarily commutative) 
$K$-algebra with some nice properties. More precisely, $A$ is 
a noetherian AS regular Koszul quiver algebra. 
If $A$ is commutative, it is essentially a polynomial ring. 
But even in this case, most results in \S4 and a few results in \S3 are new. 
(In the polynomial ring case, 
many results in \S3 were obtained in \cite{EFS}.) 
So one can read this paper assuming that $A$ is a polynomial ring. 

\medskip

We sketch the definition and basic properties of 
graded quiver algebras here. See \cite{GM} for further information.  

Let $Q$ be a finite quiver. That is, $Q = (Q_0, Q_1)$ is an oriented 
graph, where $Q_0$ is the set of vertices and $Q_1$ is the set of 
arrows. Here $Q_0$ and $Q_1$ are finite sets. 
The path algebra $KQ$ is a 
positively graded algebra with grading given by the lengths of paths. 
We denote the graded Jacobson radical of $KQ$ by $J$. That is, 
$J$ is the ideal generated by all arrows. 
If $I \subset J^2$ is a graded ideal, we say $A = KQ/I$ is a 
{\it graded quiver algebra}. Of course, $A = \bigoplus_{i \geq 0} A_i$ 
is a graded ring such that the degree $i$ component 
$A_i$ is a finite dimensional $K$-vector space for all $i$. 
The subalgebra $A_0$ is a product of copies of the field $K$, 
one copy for each element of $Q_0$.  
If $A_0 = K$ (i.e., $Q$ has only one vertex), we say $A$ is {\it connected}. 
Let  $R = \bigoplus_{i \geq 0} R_i$ be a graded algebra with 
$R_0 = K$ and $\dim_K R_1 =:n < \infty$. 
If $R$ is generated by $R_1$ as a $K$-algebra, 
then it can be regarded as a graded 
quiver algebra over a quiver with one vertex and $n$ loops. 
Let $\m := \bigoplus_{i \geq 1} A_i$ be the graded Jacobson radical of $A$.  
Unless otherwise specified, we assume that $A$ is left and right noetherian 
throughout this paper.

Let $\GrA$ (resp. $\GrAo$) be the category of graded left 
(resp. right) $A$-modules and their degree preserving $A$-homomorphisms. 
Note that the degree $i$ component $M_i$ of $M \in \GrA$ (or $M \in \GrAo$) 
is an $A_0$-module for each $i$.  
Let $\grA$ (resp. $\grAo$) be the full subcategory of $\GrA$ 
(resp. $\GrAo$) consisting of finitely generated modules. 
Since we  assume that $A$ is noetherian, $\grA$ and $\grAo$ are 
abelian categories. 
In the sequel, we will define several concepts for $\GrA$ and $\grA$. 
But the corresponding concepts for $\GrAo$ and $\grAo$ can be defined 
in the same way. 

For $n \in \ZZ$ and $M \in \GrA$, 
set $M_{\geq n} := \bigoplus_{i \geq n} M_i$ to be a 
submodule of $M$, and $M_{\leq n} := \bigoplus_{i \leq n} M_i$ 
to be a graded $K$-vector space. The $n^{\rm th}$ shift 
$M(n)$ of $M$ is defined by $M(n)_i = M_{n+i}$. 
Set $\sigma(M) := \sup \{ \, i \mid M_i \ne 0 \, \}$ and 
$\iota(M) := \inf \{ \, i \mid M_i \ne 0 \, \}$. 
If $M = 0$, we set $\sigma(M) = - \infty$ and $\iota(M) = + \infty$. 
Note that if $M \in \grA$ then $\iota(M)>-\infty$. 
For a complex $\M$ in $\GrA$, set 
$$\sigma(\M) := \sup \{ \, \sigma(H^i(\M)) + i \mid i \in \ZZ \, \}
\ \text{and} \ 
\iota(\M) := \inf \{ \, \iota(H^i(\M)) + i \mid i \in \ZZ \, \}.$$  

For $v \in Q_0$, we have the idempotent $e_v$ associated with $v$. 
Note that $1 = \sum_{v \in Q_0} e_v$. Set $P_v := A e_v $ and 
${}_v P :=  e_v A$. Then we have ${}_A A = \bigoplus_{v \in Q_0} P_v$ and 
$A_A = \bigoplus_{v \in Q_0} ({}_vP)$. 
Each $P_v$ and ${}_v P$ are indecomposable projectives. 
Conversely, any indecomposable projective in $\GrA$ (resp. $\GrAo$)
is isomorphic to  $P_v$ (resp. ${}_v P$) for some $v \in Q_0$ 
up to degree shifting. 
Set $K_v := P_v/ (\m P_v )$ and ${}_v K := {}_v P / ({}_vP \, \m)$. 
Each $K_v$ and ${}_v K$ are simple. 
Conversely, any simple object in $\GrA$ (resp. $\GrAo$) is isomorphic to 
$K_v$ (resp. ${}_vK$) for some $v \in Q_0$ up to degree shifting.  

We say a graded left (or right) $A$-module $M$ 
is {\it locally finite} 
if $\dim_K M_i < \infty$ for all $i$. If $M \in \grA$, then it is 
locally finite. 
Let $\lfA$ (resp. $\lfAo$) be the full subcategory of 
$\GrA$ (resp. $\GrAo$) consisting of locally finite modules. 

Let $C^b(\GrA)$ be the category of bounded cochain complexes in 
$\GrA$, and $D^b(\GrA)$ its derived category. 
We have similar categories for $\GrAo$, $\lfA$, $\lfAo$, $\grA$ and $\grAo$. 
For a complex $\M$ and an integer $p$, let $\M[p]$ be the $p^{\rm th}$ 
translation of $\M$. That is, $\M[p]$ is a complex with $M^i[p] = M^{i+p}$. 
Since $D^b(\grA) \cong D^b_{\grA}(\GrA) \cong D^b_{\grA}(\lfA)$, 
we freely identify these categories. 
A module $M$ can be regarded as a complex 
$\cdots \to 0 \to M \to 0 \to \cdots $ 
with $M$ at the $0^{\rm th}$ term. We can regard $\GrA$ as a full 
subcategory of $C^b(\GrA)$ and $D^b(\GrA)$ in this way.

For $M,N \in \GrA$, set $\uHom_A(M,N):= \bigoplus_{i \in \ZZ} 
\Hom_{\GrA}(M,N(i))$ to be a graded $K$-vector space with 
$\uHom_A(M,N)_i = \Hom_{\GrA}(M,N(i))$. Similarly, we can also define 
$\uHom_A^\bullet(\M,\N)$, $\bR\uHom_A(\M,\N)$, and $\uExt^i_A(\M,\N)$ 
for $\M, \N \in D^b(\GrA)$. 

If $V$ is a $K$-vector space, $V^*$ denotes the dual vector 
space $\Hom_K(V,K)$. For $M \in \GrA$ (resp. $M \in \GrAo$), 
$M^\vee := \bigoplus_{i \in \ZZ} (M_i)^*$ 
has a graded {\it right} (resp. {\it left}) $A$-module 
structure given by $(fa)(x)= f(ax)$ (resp. $(af)(x)= f(xa)$) and 
$(M^\vee)_i = (M_{-i})^*$. 
If $M \in \lfA$, then $M^\vee \in \lfAo$ and $M^{\vee \vee} \cong M$. 
In other words, $(-)^\vee$ gives exact duality functors between $\lfA$ and 
$\lfAo$, which can be extended to duality functors between 
$C^b(\lfA)$ and $C^b(\lfAo)$, or 
between $D^b(\lfA)$ and $D^b(\lfAo)$. 
In this paper, when we say $W$ is an $A$-$A$ bimodule, 
we always assume that $(aw)a' = a(wa')$ for all 
$w \in W$ and $a, a' \in A$. If $W$ is a graded $A$-$A$ bimodule, 
then so is $W^\vee$. 

It is easy to see that  $I_v :=(_vP)^\vee$ (resp. 
${}_vI := (P_v)^\vee$) is injective in $\GrA$ (resp. $\GrAo$).  
Moreover, $I_v$ and ${}_v I$ are graded injective hulls of $K_v$ and 
${}_v K$ respectively.  
In particular, the $A$-$A$ bimodule 
$A^\vee$ is injective both in $\GrA$ and in $\GrAo$. 

Let $W$ be a graded $A$-$A$-bimodule. 
For $M \in \GrA$, we can regard $\uHom_A(M,W)$ 
as a graded {\it right} $A$-module by $(fa)(x)= f(x)a$. 
We can also define $\bR\uHom_A(\M,W) \in D^b(\GrAo)$ and 
$\uExt_A^i(\M,W) \in \GrAo$ for $\M \in D^b(\GrA)$ in this way. 
Similarly,  for $\M \in D^b(\GrAo)$, 
we can make  $\bR\uHom_{\Ao}(\M,W)$ and 
$\uExt_{\Ao}^i(\M,W)$  (bounded complex of) graded left $A$-modules. 
For $M \in \GrA$, we can regard $\uHom_A(W,M)$ 
as a graded {\it left} $A$-module by $(af)(x)= f(xa)$.

For the functor $\uHom_A(-,W)$  , 
we mainly consider the case when $W =A$ or $W = A^\vee$. 
But, we have $\uHom_A(-,A^\vee) \cong (-)^\vee$. To see this, note that 
\begin{eqnarray*}
(M^\vee)_i = \Hom_K(M_{-i},K) &=& \bigoplus_{v \in Q_0} \Hom_K(e_vM_{-i},K) \\
&\cong& \bigoplus_{v \in Q_0} \Hom_K (e_vM_{-i},K_v)\\
&\cong& \Hom_{A_0} (M_{-i}, A_0).  
\end{eqnarray*}
Via the identification $(A^\vee)_0 \cong (A_0)^* \cong A_0$, 
$f \in (M^\vee)_i \cong \Hom_{A_0}(M_{-i}, A_0)$ gives a morphism 
$f': M_{\geq -i} \to A^\vee(i)$ in $\GrA$. Since 
$\Hom_\GrA (M/M_{\geq -i}, A^\vee(i)) = 0$ and $A^\vee$ is injective, 
the short exact sequence $0 \to M_{\geq -i} \to M \to M/M_{\geq -i} \to 0$ 
induces a unique extension $f'' : M \to A^\vee(i)$ of $f'$. 
From this correspondence, we have $\uHom_A(M,A^\vee) \cong M^\vee$.

Let $\P$ be a right bounded complex in $\grA$ such that each $P^i$ 
is projective. We say $\P$ is {\it minimal} 
if $d(P^i) \subset \m P^{i+1}$ for all $i$. Here $d$ is the 
differential map. Any complex $\M \in C^b(\grA)$ 
has a minimal projective resolution, that is, 
we have a minimal complex $\P$ of 
projective objects and a graded quasi-isomorphism $\P \to \M$. 
A minimal projective resolution of $\M$ is unique up to isomorphism.
We denote a graded module $A/\m$ by $A_0$.  
Set $\beta^{i,j}(\M) := \dim_K \uExt_A^{-i}(\M,A_0)_{-j}$. 
Let $\P$ be a minimal projective resolution of $\M$, 
and $P^i := \bigoplus_{l=1}^m T^{i, \, l}$ an indecomposable decomposition. 
Then we have 
$$\beta^{i,j}(\M) = \# \{ \, l \mid \text{$T^{i, \, l}(j) 
\cong P_v$ for some $v$}\, \}.$$
We can also define $\beta^{i,j}(\M)$ as the dimension of 
$\Tor^A_{-i}(A_0, \M)_j$. This definition must be much more familiar to 
commutative algebraists. Note that $\beta^{i,j}(-)$ is an invariant of 
isomorphism classes of the derived category $D^b(\grA)$. Note that 
these facts on minimal projective resolutions also hold over any 
noetherian graded algebra. 

\begin{dfn}\label{AS regular} 
Let $A$ be a (not necessarily noetherian) graded quiver algebra. 
We say $A$ is {\it Artin-Schelter regular} 
(AS-regular, for short), if 
\begin{itemize}
\item $A$ has finite global dimension $d$. 
\item $\uExt_A^i(K_v, A) = \uExt_{\Ao}^i({}_v K,A)= 0$ for all $i \ne d$ 
and all $v \in Q_0$.  
\item  There are a permutation $\delta$ on $Q_0$ and 
an integer $n_v$ for each $v \in Q_0$ such that  
$\uExt_A^d(K_v, A) \cong {}_{\delta(v)} K(n_v)$ (equivalently, 
$\uExt_{\Ao}^d({}_v K, A) \cong K_{\delta^{-1}(v)}(n_v)$ ) 
for all $v$.   
\end{itemize}
\end{dfn}

\begin{rem}\label{original definition}
The AS regularity is a very important concept in 
non-commutative algebraic geometry. 
In the original definition,  it is assumed that 
an AS regular algebra $A$ is connected and  
there is a positive real number $\gamma$ such that 
$\dim_K A_n < n^\gamma$ for $n \gg 0$, while some authors do not require 
the latter condition. We also remark that 
Martinez-Villa and coworkers called rings satisfying 
the conditions of Definition~\ref{AS regular} {\it generalized Auslander 
regular algebras} in \cite{GMT, M98}. 
\end{rem}

\begin{dfn}
For an integer $l \in \ZZ$, 
we say $\M \in \grA$ has an {\it $l$-linear (projective) resolution}, 
if $$\beta^{i,j}(\M) \ne 0 \Rightarrow i+j=l.$$
If $\M$ has an $l$-linear resolution for some $l$, 
we say $\M$ has a {\it linear resolution}. 
\end{dfn}

\begin{dfn}
We say $A$ is {\it Koszul}, if the graded left 
$A$-module $A_0$ has a linear resolution. 
\end{dfn}

In the definition of the Koszul property, we can regard $A_0$ as a 
right $A$-module. (We get the equivalent definition.)  
That is, $A$ is Koszul if and only if any simple   
graded left (or, right) $A$-module has a linear resolution. 

\begin{lem}
If $A$ is noetherian, AS-regular, Koszul, and has global dimension $d$, 
then $\uExt_A^d(K_v, A) \cong {}_{\delta(v)} K (d)$ and 
$\uExt_{\Ao}^d({}_v K, A) \cong K_{\delta^{-1} (v)}(d)$ for all $v$.  
Here $\delta$ is the permutation of $Q_0$ given in 
Definition~\ref{AS regular}. 
\end{lem}

\begin{proof}
Since $A$ is Koszul, $P^{-d}$ of 
a minimal projective resolution $\P : 0 \to P^{-d} \to \cdots \to P^0 \to 0$ 
of $K_v$ is generated by its degree $d$-part $(P^{-d})_d$ 
(more precisely, $P^{-d} = P_{\delta(v)}(-d)$). 
\end{proof}

{\it In the rest of this paper, $A$ is always a noetherian AS-regular 
Koszul quiver algebra of global dimension $d$.}

\begin{exmp}
(1) A polynomial ring $K[x_, \ldots, x_d]$ is clearly 
a noetherian AS-regular Koszul (quiver) algebra of global dimension $d$.
Conversely, if a regular noetherian graded algebra 
is connected and commutative, it is a polynomial ring. 

(2) Let $K\< x_1, \ldots, x_d\>$ be the free associative algebra, and 
$(q_{i,j})$ a $d \times d$ matrix with entries in $K$ satisfying 
$q_{i,j} q_{j,i} = q_{i,i} = 1$ for all $i, j$. Then the quotient ring 
$A = K\< x_1, \ldots, x_n\>/\<x_j x_i -q_{i,j} x_i x_j \mid 
1 \leq i, j \leq d\>$ is a noetherian AS-regular Koszul algebra with 
global dimension $d$.  This fact must be well-known to specialists, 
but we will sketch a proof here for the reader's convenience. 
Since $x_1, ..., x_d \in A_1$ form a regular 
normalizing sequence with the quotient ring $K = A/(x_1, \ldots, x_d)$, 
$A$ is a noetherian ring with a balanced dualizing complex 
by \cite[Lemma 7.3]{Mori}.  
It is not difficult to construct a minimal free resolution of the module 
$K = A/\m$, which is a ``$q$-analog" of the Koszul complex of 
a polynomial ring $K[x_1, \ldots, x_d]$. So $A$ is Koszul 
and has global dimension $d$.  
Since $A$ has finite global dimension and admits 
a balanced dualizing complex,  
it is AS-regular (c.f. \cite[Remark 3.6 (3)]{Mori}).

Artin, Tate and Van den Bergh (\cite{ATV}) classified connected 
AS regular algebras of global dimension 3.  
(Their definition of AS regularity is stronger than ours. 
See Remark~\ref{original definition}.) All of the algebras they listed   
are noetherian (\cite[Theroem~8.1]{ATV}). 
But some are Koszul and some are not.

(3) A {\it preprojective algebra} is an important example of 
non-connected AS regular algebras. 
See \cite{GMT} and the references cited there for the definition 
of this algebra and further information. 
The preprojective algebra $A$ of a finite quiver $Q$ is a graded quiver 
algebra over the {\it inverse completion} $\overline{Q}$ of $Q$. 
If the quiver $Q$ is connected (of course, it does not mean $A$ is connected), 
then $A$ is (almost) always an AS regular algebra of global dimension 2, 
but it is not Koszul in some cases, and not noetherian in many cases. 
Let $G$ be the bipartite graph of $Q$ in the sense of 
\cite[\S3]{GMT}. If $G$ is Euclidean, then $A$ is a noetherian 
AS-regular Koszul algebra of global dimension 2.  
\end{exmp}

 For $M \in \GrA$, set 
 $$\Gamma_\m(M) = \nlim \uHom_A(A/\m^n, M) 
 = \{ \, x \in M \mid \text{$A_{n} \, x =0$ for $n \gg 0$} \, \}\in \GrA.$$  
Then $\Gamma_\m(-)$ gives a left exact functor from $\GrA$ to itself.
So we have a right derived functor $\bR \Gamma_\m : D^b(\GrA) \to 
D^b(\GrA)$.  For $\M \in D^b(\GrA)$, $H_\m^i(\M)$ denotes the $i^{\rm th}$ 
cohomology of $\bR \Gamma_\m(\M)$, and we call it the $i^{\rm th}$  
{\it local cohomology} of $\M$. 
It is easy to see that $H_\m^i(\M) =  \nlim \uExt^i_A(A/\m^n, \M)$. 
Similarly, we can define $\bR \Gamma_{\m^\op}: D^b(\GrAo) \to 
D^b(\GrAo)$ and $H^i_{\m^\op}: D^b(\GrAo) \to \GrAo$ in the same way. 
If $M$ is an $A$-$A$ bimodule, $H_\m^i(M)$ and $H_{\m^\op}^i(M)$ 
are also. 

Let $I \in \GrA$ be an indecomposable injective.  
Then $\Gamma_\m(I) \ne 0$, if and only if $I \cong I_v(n)$ 
for some $v \in Q_0$ and $n \in \ZZ$, if and only if 
$\Gamma_\m(I) = I$.  Similarly, $\uHom_A(A_0, I) \ne 0$ if and 
only if  $I \cong I_v(n)$ for some $v \in Q_0$ and $n \in \ZZ$. In this case, 
$\uHom_A(A_0, I) = K_v(n)$.  
The same is true for an indecomposable injective $I \in \GrAo$.

Let $\I$ be a minimal injective resolution of $A$ in $\grA$. 
Since $A$ is AS regular, $I^i = 0$ for all $i>d$, $\Gamma_\m(I^i) = 0$ 
for all $i < d$, and $\Gamma_\m(I^d) = A^\vee(d)$. 
Hence $\bR\Gamma_\m(A) \cong A^\vee(d)[-d]$ in $D^b(\grA)$. 
By the same argument as \cite[Proposition~4.4]{Ye}, we also have 
$\bR\Gamma_\m(A) \cong A^\vee(d)[-d]$ in $D^b(\grAo)$. It does not mean that 
$H_\m^d(A) \cong A^\vee(d)$ as $A$-$A$ bimodules. 
But there is an $A$-$A$ bimodule $L$ such that 
$L \otimes_A H_\m^d(A)  \cong A^\vee(d)$ as $A$-$A$ bimodules. 
Here the underlying additive group of $L$ is $A$, 
but the bimodule structure is give by 
$A \times L \times A \ni (a, l,  b) \mapsto \phi(a) l b \in A = L$ 
for a (fixed) $K$-algebra automorphism $\phi$ of $A$. 
In particular, $L \cong A$ as 
left $A$-modules and as right $A$-modules (separately). 
Note that $\phi(e_v) = e_{\delta(v)}$ for all $v \in Q_0$, where 
$\delta$ is the permutation on $Q_0$ appeared in 
Definition~\ref{AS regular}. If $A$ is commutative, then 
$\phi$ is the identity map. 

We give a new $A$-$A$ bimodule structure $L'$ to the additive group $A$ by  
$A \times L' \times A \ni (a, l,  b) \mapsto a l \phi(b) \in A = L'$. 
Then $L' \cong \uHom_A(L,A)$.  Set $\D := L'(-d)[d]$. 
Note that $\D$ belongs both $D^b(\grA)$ and $D^b(\grAo)$.  We have 
$H_\m^i(\D) = H_{\m^\op}^i(\D) = 0$ for all $i \ne 0$ and 
$H^0_\m (\D) \cong H_{m^\op}^0(\D) \cong A^\vee$ as $A$-$A$ 
bimodules by the same argument as \cite[\S4]{Ye}. 
Thus (an injective resolution of) $\D$ is a 
{\it balanced dualizing complex} of $A$ in the sense of \cite{Ye} 
(the paper only concerns connect rings, but the definition 
can be generalized in the obvious way).

Easy computation shows that 
$\uHom_A(P_v, L') \cong {}_{\delta^{-1}(v)} P$ and $\uHom_{\Ao}
( {}_vP, L') \cong P_{\delta(v)}$ for all $v \in Q_0$. Since   
$\bR \uHom_A(\M, \D)$ (resp. $\bR \uHom_{\Ao}(\M, \D)$) for 
$\M \in \grA$ (resp. $\M \in \grAo$) can be computed by a projective 
resolution of $\M$, $\bR \uHom_A(-, \D)$ and $\bR \uHom_{\Ao}(-, \D)$ 
give duality functors between $D^b(\grA)$ and $D^b(\grAo)$.  
(Of course, we can also prove this by 
the same argument as \cite[Proposition~3.4]{Ye}.)

\begin{thm}[{Yekutieli \cite[Theorem 4.18]{Ye}, Martinez-Villa 
\cite[Proposition~4.6]{M98}}]\label{Serre duality}
For $\M \in D^b(\grA)$, we have 
$$\bR \Gamma_\m(\M)^\vee \cong \bR\uHom_A(\M, \D).$$
In particular, $$(H^i_\m(\M)_j)^* \cong \uExt_A^{-i}(\M, \D)_{-j}.$$
\end{thm}

\begin{proof}
The above result was proved by Yekutieli in the connected case. 
(In some sense, Martinez-Villa proved a more general result than ours, 
but he did not concern complexes.) But, the proof of  
\cite[Theorem 4.18]{Ye} only uses formal properties such as  
$A$ is noetherian, $\bR\uHom_{\Ao} (\bR\uHom_A(-, \D), \, \D) \cong \Id$, 
and $\bR\Gamma_\m \D \cong A^\vee$.  
So the proof also works in our case. 
\end{proof}

\begin{dfn}[J\o rgensen, \cite{Jor04}]
For $\M \in D^b(\grA)$, we say  
$$\reg(\M) := \sigma(\bR \Gamma_\m(\M))
= \sup \{ \, i+j \mid H_\m^i(\M)_j \ne 0 \, \}$$ 
is the {\it Castelnuovo-Mumford regularity} of $\M$. 
\end{dfn}

By Theorem~\ref{Serre duality} and the fact that 
$\bR \uHom_A(\M, \D) \in D^b(\grAo)$, 
we have $\reg(\M) < \infty$ for all $\M \in D^b(\grA)$. 

\begin{thm}[J\o rgensen, \cite{Jor04}]\label{CM reg = Ext reg} 
If $\M \in C^b(\grA)$, then 
\begin{equation}\label{reg}
\reg(\M) = \max \{ \, i+j \mid \beta^{i,j}(\M) \ne 0 \, \}.
\end{equation}
\end{thm}

When $A$ is a polynomial ring and $\M$ is a module, the above theorem 
is a fundamental result obtained by Eisenbud and Goto \cite{EG}. 
In the non-commutative case, 
under the assumption that $A$ is connected but not necessarily regular, 
this has been proved by J\o rgensen \cite[Corollary~2.8]{Jor04}. 
(If $A$ is not regular, we have $\reg(A) > 0$ in many cases.  
So one has to assume that $\reg A = 0$ there.) 
In our case (i.e., $A$ is AS-regular), we have a much simpler 
proof. So we will give it here. This proof is also different from one given in 
\cite{EG}. 

\begin{proof}
Set $\Q := \uHom_A^\bullet(\P, L'(-d)[d] \, )$. Here 
$\P$ is a minimal projective resolution of $\M$, and $L'$ is the $A$-$A$ 
bimodule given in the construction of the dualizing complex $\D$. 
Recall that $\uHom_A(P_v, L') \cong {}_{\delta^{-1}(v)} P$ for all 
$v \in Q_0$. Let $s$ be the right hand side of \eqref{reg}, and 
$l$ the minimal integer 
with the property that $\beta^{l, s-l}(\M) \ne 0$. 
Then $\iota(Q^{-d-l}) = l-s+d$, and $(Q^{-d-l+1})_{\leq (l-s+d-1)} = 0$ 
(Note that $\beta^{l-1,m}(\M) = 0$ for all $m \geq s -l +1$). 
Since $\Q$ is a minimal complex, we have 
$$0 \ne H^{-d-l}(\Q)_{l-s+d} = \uExt_A^{-d-l}(\M,\D)_{l-s+d} 
= (H^{d+l}_{\m}(\M)_{-l+s-d})^*.$$ 
Thus $\reg(\M) \geq \max \{ \, i+j \mid \beta^{i,j}(\M) \ne 0 \, \}$. 

On the other hand,  if $H^{d+l}_{\m}(\M)_{-l+r-d} \ne 0$, 
we have that $\beta^{l, t-l}(\M) \ne 0$ for some $t \geq r$ by an argument 
similar to the above. Hence $\reg(\M) \leq 
\max \{ \, i+j \mid \beta^{i,j}(\M) \ne 0 \, \}$, and we are done. 
\end{proof}

For $\M \in D^b(\grA)$, set $\cH(\M)$ to be the complex 
such that $\cH(\M)^i = H^i(M)$ for all $i$ 
and all differential maps are zero. 

\begin{lem}\label{degenerate} 
We have $\beta^{i,j}(\cH(\M)) \geq \beta^{i,j}(\M)$ 
for all $\M \in D^b(\grA)$ and 
all $i,j \in \ZZ$. In particular, $\reg(\cH(\M)) \geq \reg(\M)$. 
\end{lem}

The difference between $\reg(\M)$ and $\reg(\cH(\M))$ can be 
arbitrary large. In the last section, we will study 
the relation between this difference and a work of Martinez-Villa 
and Zacharia~\cite{MZ}.  

\begin{proof}
The assertion easily follows from the spectral sequence 
$$E_2^{p,q} = \uExt_A^p(H^{-q}(\N),A_0) \too \uExt^{p+q}_A(\N,A_0).$$
\end{proof}

For a complex $\M \in C^b(\grA)$ and an integer $r$, $(\M)_{\geq r}$ 
denotes the subcomplex of $\M$ whose $i^{\rm th}$ term is 
$(M^i)_{\geq (r-i)}$. 
Even if $\M \cong \N$ in $D^b(\grA)$, we have $(\M)_{\geq r} \not \cong 
(\N)_{\geq r}$ in general. 

In the module case, the following is a well-known property of 
Castelnuovo-Mumford regularity. 

\begin{prop}\label{truncate}
Let $\M \in C^b(\grA)$. Then $(\M)_{\geq r}$ has an 
$r$-linear resolution if and only if $r \geq \reg(\M)$.   
\end{prop}

To prove the proposition, we need the following lemma. 

\begin{lem}\label{artinian}
For a module $M \in \grA$ with $\dim_K M < \infty$, we have 
$H^0_\m(M)=M$ and $H^i_\m(M)=0$ for all $i \ne 0$. 
In particular, $\reg(M) = \sigma(M)$ in this case. 
\end{lem} 

\begin{proof}
If $\P$ is a minimal projective resolution of $M^\vee \in \grAo$, then  
$\I := (\P)^\vee$ is a minimal injective resolution of $M$. 
Since each indecomposable summand of $I^i$ is isomorphic to $I_v(n)$ for 
some $v \in Q_0$ and $n \in \ZZ$, we have $\Gamma_\m(\I) = \I$.   
\end{proof}

\noindent{\it Proof of Proposition~\ref{truncate}.}
For a complex $\L \in D^b(\grA)$, it is easy to see that 
$\iota(\L) = \min \{ \, i+j \mid \beta^{i,j}(\L) \ne 0 \, \}$. 
In particular, $\iota(\L) \leq \reg(\L)$. Hence $\L$ has an $l$-linear 
projective resolution if and only if $\iota(\L) = \reg(\L) = l$. 

Consider the short exact sequence of complexes 
\begin{equation}\label{find the way}
0 \to (\M)_{\geq r} \to \M \to \M/(\M)_{\geq r} \to 0, 
\end{equation} and set $\N:= \M/(\M)_{\geq r}$.  Note that 
$\dim_K H^i(N) < \infty$ for all $i$.  
By Lemmas~\ref{artinian} and \ref{degenerate}, 
we have $$r > \sigma(\N) = \max\{ \, \reg(H^i(\N))+i \mid i \in \ZZ  \, \} 
= \reg(\cH(\N)) \geq \reg(\N).$$ 
By the long exact sequence of $\uExt_A^\bullet(-, A_0)$ induced by 
\eqref{find the way}, we have 
\begin{eqnarray*}
r \leq \iota( (\M)_{\geq r} ) 
\leq \reg((\M)_{\geq r}) 
&\leq&  \max \{ \, \reg(\N)+1, \, \reg(\M) \, \} \\ 
&\leq& \max \{ \, r, \, \reg(\M) \, \}.
\end{eqnarray*} 
Moreover, if $r < \reg(\M)$ then we have $\reg (\N)+1 < \reg(\M)$ and 
$\reg((\M)_{\geq r}) = \reg(\M) > r$. 
Hence $(\M)_{\geq r}$ has an $r$-linear resolution if and only if 
$r \geq \reg(\M)$.  
\qed

\medskip

The following is one of the most basic 
results on Castelnuovo-Mumford regularity (see \cite{EG}). 
J\o rgensen~\cite{Jor99} proved the same result for $M \in \grA$. 

\medskip

{\it Let $S = K[x_1, \ldots, x_d]$ be a polynomial ring. 
If $M \in \grS$ satisfies $H_\m^0(M)_{\geq r+1} = 0$ and 
$H_\m^i(M)_{r+1-i}=0$ for all $i \geq 1$, then $r \geq \reg(M)$ 
(i.e., $H_\m^i(M)_{\geq r+1-i}=0$ for all $i \geq 1$). }

\medskip

The similar result also holds for $\M \in D^b(\grA)$. 
Since a minor adaptation of the proof of 
\cite[Theorem~2.4]{Jor99} also works for complexes, 
we leave the proof to the reader. 

\begin{prop} 
If $\M \in D^b(\grA)$ with 
$t := \max\{ \, i \mid H^i(\M) \ne 0  \, \}$ satisfies 
\begin{itemize}
\item $H_\m^i(\M)_{\geq r+1-i} = 0$ for all $i \leq t$ 
\item $H_\m^i(\M)_{r+1-i} = 0$ for all $i > t$, 
\end{itemize}
then $r \geq \reg(\M)$ (i.e., $H_\m^i(\M)_{\geq r+1-i} = 0$ for all $i > t$).
\end{prop}

\section{Koszul duality}
In this section, we study the relation between the Castelnuovo-Mumford 
regularity of complexes and the Koszul duality. 
For precise information of this duality, see \cite[\S2]{BGS}.
There, the symbol $A$ (resp. $A^!$) basically 
means a finite dimensional (resp. noetherian) Koszul algebra. 
This convention is opposite to ours. So the reader should be careful. 

\medskip

Recall that $A = KQ/I$ is a graded quiver algebra over a finite quiver $Q$. 
Let $Q^\op$ be the {\it opposite quiver} of $Q$. That is, $Q_0^\op = Q_0$ and 
there is a bijection from $Q_1$ to $Q^\op_1$ which sends an arrow 
$\alpha: v \to u$ in $Q_1$ to the arrow 
$\alpha^\op : u \to v$ in $Q_1^\op$. 
Consider the bilinear form $\< -, - \> : (KQ)_2 \times (KQ^\op)_2 
\to K$ defined by 
$$
\< \alpha \beta, \gamma^\op \delta^\op \> = \begin{cases}
1 & \text{if $\alpha = \delta$ and $\beta = \gamma$,} \\
0 & \text{otherwise}
\end{cases}
$$
for all $\alpha, \beta, \gamma, \delta \in Q_1$. 
Let $I^\bot \subset KQ^\op$ be the ideal generated by 
$$\{ \, y \in (KQ^\op)_2 \mid 
\text{$\<x, y\> = 0$ for all $x \in I_2$} \, \}.$$ 
We say $KQ^\op/I^\bot$ is the  {\it quadratic dual ring} of $A$, 
and denote it by $A^!$.  Clearly, $(A^!)_0 = A_0$.  
Since $A$ is Koszul, so is $A^!$ and 
$(A^!)^! \cong A$. Since $A$ is AS regular, $A^!$ is a finite dimensional 
selfinjective algebra with $A = \bigoplus_{i=0}^d A_i$ 
by \cite[Theorem~5.1]{M99}. If $A$ is a polynomial ring, then 
$A^!$ is the exterior algebra $\bigwedge (A_1)^*$.  

Since $A^!$ is selfinjective, 
$\DAq := \uHom_{A^!}(-, A^!)$ and $\DAqo := \uHom_{\Aqo}(-, A^!)$
give exact duality functors between $\grAq$ and $\grAqo$. 
They induce duality functors between  
$D^b(\grAq)$ and $D^b(\grAqo)$, which are also denoted by $\DAq$and $\DAqo$. 
It is easy to see that $\DAq(N) \cong \Hom_K(N,K)(-d)$.

We say a complex $\F \in C(\grAq)$ is a projective 
(resp. injective) resolution of a complex $\N \in C^b(\grAq)$, 
if each term $F^i$ is projective ($=$ injective), 
$\F$ is right (resp. left) bounded, and there is a graded quasi isomorphism 
$\F \to \N$ (resp. $\N \to \F$). 
We say a projective (or, injective) resolution $\F \in C^b(\grAq)$ is 
{\it minimal} if $d^i(F^i) \subset \n F^{i+1}$ for all $i$, 
where $\n$ is the graded Jacobson radical of $A^!$. 
(The usual definition of a minimal injective resolution is different 
from the above one. But they coincide in our case.) 
A bounded complex $\N \in C^b(\grAq)$ has a minimal projective 
resolution and a minimal injective resolution, and they are unique 
up to isomorphism. If $\F$ is a minimal projective (resp. injective) 
resolution of $\N$ then $\DAq(\F)$ is a  minimal injective (resp. projective) 
resolution of $\DAq(\N)$.  

For $\N \in D^b(\grAq)$, set 
$$\mu^{i,j}(\N):= \dim_K \uExt_{A^!}^i(A_0, \N)_j.$$ 
Then $\mu^{i,j}(\N)$  measures the size of 
a minimal injective resolution of $\N$. More precisely, 
if $\F$ is a minimal injective resolution of $\N$, 
and $F^i := \bigoplus_{l =1}^m T^{i,\, l}$ is an indecomposable 
decomposition, then we have 
\begin{eqnarray*}
\mu^{i,j}(\N) &=& \# \{ \, l \mid \text{
$\operatorname{soc}(T^{i,\, l}) = (T^{i,\, l})_j$ }\, \}\\
&=& \# \{ \, l \mid  \, \text{$T^{i,\, l}(j)$ 
is isomorphic to a direct summand of $A^!(d)$} \, \}.
\end{eqnarray*} 

\medskip

Let $V$ be a finitely generated left $A_0$-module. Then  
$\Hom_{A_0}(A^!, V)$ is a graded  {\it left}  $A^!$-module with 
$(af)(a') = f(a'a)$ and $\Hom_{A_0}(A^!, V)_i = 
\Hom_{A_0} ((A^!)_{-i}, V)$. 
Since $A^!$ is selfinjective, we have 
$\Hom_{A_0}(A^!, A_0) \cong A^!(d)$. Hence $\Hom_{A_0}(A^!,V)$ is 
a projective (and injective) left $A^!$-module for all $V$. 
If $V$ has degree $i$ (e.g.,  $V = M_i$ for some $M \in \grA$), then 
we set  $\Hom_{A_0}(A^!, V)_j = \Hom_{A_0}(A^!_{-j-i}, V)$. 

For $\M \in C^b(\grA)$, let 
$\DG(\M):= \Hom_{A_0}(A^!,\M) \in C^b(\grAq)$ be  
the total complex of the double complex with 
$\DG(\M)^{i,j} =  \Hom_{A_0}(A^!, M^i_j)$ whose vertical and horizontal 
differentials $d'$ and $d''$ are defined by 
$$d'(f)(x)= \sum_{\alpha \in Q_1} \alpha f(\alpha^\op x),  
\qquad \quad d''(f)(x) = \partial_{\M}(f(x))$$
for $f \in \Hom_{A_0}(A^!, M^i_j)$ and $x \in A^!$.  
The gradings of $\DG(\M)$ is given by 
$$\DG(\M)^p_q := \bigoplus_{p=i+j, \, q = -l-j} 
\Hom_{A_0}((A^!)_l, M^i_j).$$ 

Each term of $\DG(\M)$ is injective. 
For a module $M \in \grA$, $\DG(M)$ is a minimal complex. 
Thus we have 
\begin{equation}
\mu^{i,j}(\DG(M)) = \begin{cases} 
\dim_K M_i & \text{if $i+j = 0$,}\\
0 & \text{otherwise.}
\end{cases}
\end{equation}

Similarly, for a complex $\N \in C^b(\grAq)$, we can define a new complex 
$\DF(\N):= A \otimes_{A_0} \N \in C^b(\grA)$ as the total complex 
of the double complex with 
$\DF(\N)^{i,j} =  A \otimes_{A_0} N^i_j$ whose vertical and horizontal 
differentials $d'$ and $d''$ are defined by 
$$d'(a \otimes x)= \sum_{\alpha \in Q_1} a \alpha \otimes  \alpha^\op x, 
\qquad d''(a \otimes x)=  a \otimes \partial_{\N}(x) $$
for $a \otimes x \in A \otimes_{A_0} N^i$. The gradings of $\DF(\N)$
is given by 
$$\DF(\N)^p_q:= \bigoplus_{p=i+j, \, q = l-j} A_l \otimes_{A_0} N^i_j.$$
Clearly, each term of $\DF(\N)$ is a projective $A$-module.  
For a module 
$N \in \grAq$,  $\DF(N)$ is a minimal complex. Hence we have 
\begin{equation}\label{betti of DF}
\beta^{i,j}(\DF(N)) = \begin{cases} 
\dim_K N_i & \text{if $i+j = 0$,}\\
0 & \text{otherwise.}
\end{cases}
\end{equation}

It is well-known that the operations $\DF$ and $\DG$ 
define functors $\DF : D^b(\grAq) \to D^b(\grA)$ and 
$\DG:  D^b(\grA) \to D^b(\grAq)$, and they give an equivalence 
$D^b(\grA) \cong D^b(\grAq)$ of triangulated categories. 
This equivalence is called the {\it Koszul duality}. 
When $A$ is a polynomial ring, 
this equivalence is called {\it Bernstein-Gel'fand-Gel'fand 
correspondence}. See, for example, \cite{EFS}. 

Note that $(A^\op)^! \cong (A^!)^\op$. 
We have the functors $\DFo : D^b(\grAqo) \to D^b(\grAo)$ 
and $\DGo: D^b(\grAo) \to D^b(\grAqo)$ giving  
$D^b(\grAo) \cong D^b(\grAqo)$.

\begin{prop}[{c.f. \cite[Proposition~2.3]{EFS}}]
\label{generalization of EFS} In the above situation, we have 
$$\beta^{i,j}(\M) = \dim_K H^{i+j}(\DG(\M))_{-j} \quad  \text{and} \quad  
\mu^{i,j}(\N) = \dim_K H^{i+j}(\DF(\N))_{-j}.$$
\end{prop}

\begin{proof} 
While the assertion follows from Proposition~\ref{strand} below, 
we give a direct proof here. We have 
\begin{eqnarray*}
\uExt_{A^!}^i(A_0, \N)_j &\cong& 
\Hom_{D^b(\grAq)}(A_0, \N[i](j))\\
&\cong& \Hom_{D^b(\grA)}( \, \DF(A_0), \, \DF(\N[i](j)) \, )\\ 
&\cong& \Hom_{D^b(\grA)}( \, A, \, \DF(\N)[i+j](-j) \, )\\
&\cong& H^{i+j}(\DF(\N))_{-j}. 
\end{eqnarray*}
Since $\mu^{i,j}(\N)= \dim_K \uExt_{A^!}^i(A_0, \N)_j$,  
the second equation of the proposition follows. 
We can prove the first equation by a similar argument. 
But in this time we use the contravariant functor 
$\DAq \circ \DG : D^b(\grA) \to D^b(\grAqo)$ 
and the fact that $\DAq \circ \DG(A_0) \cong 
\DAq(A^!(d)) \cong A^!(-d)$. 
\end{proof}

\begin{cor}\label{reg & DG}
$\reg(\M) = \max \{ \, i \mid  H^i(\DG(\M)) \ne 0 \, \}$.
\end{cor}

\begin{proof}
Follows Theorem~\ref{CM reg = Ext reg} 
and Proposition~\ref{generalization of EFS}.
\end{proof}

Recall that $\DA := \bR \uHom_A(-,\D)$ is a duality functor 
from $D^b(\grA)$ to $D^b(\grAo)$. 

\begin{prop}\label{reg of dual}
$\reg( \, \DA(\M) \, ) = - \min \{ \, i \mid  H^i(\DG(\M)) \ne 0 \, \}$. 
\end{prop}

\begin{proof}
Let $L'$ be the $A$-$A$ 
bimodule given in the construction of the dualizing complex $\D$. 
Note that $\DA(\M) \cong \uHom_A^\bullet(\P,L'(-d)[d]) =:\Q$ 
for a projective resolution $\P$ of $\M$. 
Since $\DA(P_v) = {}_{\delta^{-1}(v)} P(-d)[d]$, $\Q$ 
is a complex of projectives. And $\Q$ 
is a minimal complex if and only if $\P$ is. Hence 
$\beta^{-i-d, -j+d}(\DA(\M)) = \beta^{i,j}(\M)$. Therefore, the assertion 
follows from Proposition~\ref{generalization of EFS}. 
\end{proof}

We can refine Proposition~\ref{generalization of EFS} using 
the notion of {\it linear strands} of projective (or injective) resolutions, 
which was introduced by Eisenbud et. al. (See \cite[\S3]{EFS}.) 
First, we will generalize this notion to our rings. 
Let $B$ be a noetherian Koszul algebra (e.g. $B= A$ or $A^!$) 
with the graded Jacobson radical $\m$, 
and $\P$ a {\it minimal} projective resolution of a bounded complex 
$\M \in D^b(\grB)$.  
Consider the decomposition $P^i := \bigoplus_{j \in \ZZ} P^{i,j}$ 
such that any indecomposable summand of $P^{i,j}$ 
is isomorphic to a summand of $B(-j)$. 
For an integer $l$, we define the {\it $l$-linear strand} 
$\plin_l(\M)$ of a projective resolution of 
$\M$ as follows: The term $\plin_l(\M)^i$ of cohomological degree $i$  
is $P^{i, l-i}$ and the differential 
$P^{i, l-i} \to P^{i+1, l-i-1}$ 
is the corresponding component of the differential 
$P^i \to P^{i+1}$ of $\P$.  So the differential of $\plin_l(\M)$ 
is  represented by a matrix whose entries are elements in $A_1$. 
Set $\plin(\M) := \bigoplus_{l \in \ZZ} \plin_l(\M)$. 
It is obvious that $\beta^{i,j}(\M) = \beta^{i,j}(\plin(\M))$ 
for all $i,j$. 

Using spectral sequence argument, we can construct $\plin (\M)$ from  
a (not necessarily minimal) projective resolution $\Q$ of $\M$. 
Consider the $\m$-adic filtration 
$\Q = F_0\Q \supset F_1 \Q \supset \cdots$ of 
$\Q$ with $F_p Q^i = \m^p Q^i$ and the associated spectral 
sequence $\{E_r^{*,*},d_r \}$. The associated graded object 
$\grm M := \bigoplus_{p \geq 0} \m^p 
M/\m^{p+1}M$ of $M \in \grB$ is a module over 
$\grm B = \bigoplus_{p \geq 0} \m^p/\m^{p+1} \cong B$. 
Since $\m^p M$ is a graded submodule of $M$, we can make $\grm M$ 
a graded module using the original grading of $M$ 
(so $(\grm M)_i$ is {\it not} $\m^i M /\m^{i+1} M$ here). 
Under the identification $\grm B$ with $B$,  
we have $\grm M \not \cong M$ in general. 
But if each indecomposable summand $N$ of $M$ is generated by 
$N_{\iota(N)}$ then $\grm M \cong M$. 
Since $Q^t$ is a projective $B$-module, $Q_0^t := \bigoplus_{p+q=t}
E_0^{p,q} = \bigoplus_{p \geq 0} \m^p Q^t / \m^{p+1} Q^t = \grm Q^t$ 
is isomorphic to $Q^t$. 
The maps $d_0^{p,q} : E_0^{p,q} \to E_0^{p,q+1}$ make $\Q_0$ 
a cochain complex of projective $\grm B$-modules.  
Consider the decomposition $\Q = \P \oplus C^\bullet$, where $\P$ is minimal 
and $C^\bullet$ is exact. (We always have such a decomposition.) 
If we identify $Q_0^t$ with 
$Q^t = P^t \oplus C^t$, the differential $d_0$ of $\Q_0$ is given by 
$(0, d_{C^\bullet})$. 
Hence we have $Q_1^t = \bigoplus_{p+q=t}E_1^{p,q} \cong P^t$. 
The maps $d_1^{p,q}: E_1^{p,q} = \m^p P^t/\m^{p+1} P^t  
\to E_1^{p+1,q} = \m^{p+1} P^{t+1}/\m^{p+2} P^{t+1}$ make 
$\Q_1$ a cochain complex of projective $\grm B (\cong B)$-modules 
whose differential is the ``linear component" of the differential 
$d_{\P}$ of $\P$. Thus the complex 
$(Q_1^\bullet, d_1)$ is isomorphic to  $\plin (\M)$. 

\medskip

Since $A^!$ is selfinjective, we can consider the linear strands 
of an injective resolution. More precisely, 
starting from a minimal injective resolution of $\N \in D^b(\grAq)$, 
we can construct its $l$-linear strand 
$\ilin_l(\N)$  in a similar way. Here, if $I^i$ 
is the cohomological degree $i^{\rm th}$ term of $\ilin_l(\N)$, 
then the socle of $I^i$ coincides with $(I^i)_{l-i}$. In other words, 
any indecomposable summand of $I^i$ 
is isomorphic to a summand of $A^!(i-l+d)$.  
Set $\ilin (\N) = \bigoplus_{l \in \ZZ} \ilin_l(\N)$. 
This complex can also be constructed using spectral sequences. 

We have that  $\DAq(\ilin(\N)) \cong \plin(\DAq(\N))$ 
and $\DAq(\plin(\N)) \cong \ilin(\DAq(\N))$. 

\begin{prop}[{c.f. \cite[Corollary~3.6]{EFS}}]\label{strand} 
For  $\M \in D^b(\grA)$ and $\N \in D^b(\grAq)$, 
we have $$\plin (\DF(\N)) = \DF(\cH(\N)) \quad  \text{and} \quad   
\ilin (\DG(\M)) = \DG(\cH(\M)).$$ More precisely, $$\plin_l (\DF(\N)) = 
\DF(H^l(\N))[-l] \quad  \text{and} \quad  
\ilin_l(\DG(\M)) = \DG(H^l(\M))[-l].$$ 
\end{prop}

\begin{proof}
Set $\Q = \DF(\N)$. Note that $\Q$ is a (non minimal) complex 
of projective modules. We use the above spectral sequence argument 
(and the notation there). Then the differential  
$d_0^t : Q_0^t \cong \DF^t(\N) \to Q_0^{t+1} \cong 
\DF^{t+1}(\N)$ is given by $\pm \partial_{\N}$. Thus  
$$Q_1^t \cong \bigoplus_{t = i+j} A \otimes_{A_0} H^i(\N)_j 
= \bigoplus_{t = i+j}\DF^j(H^i(\N)),$$ 
and the differential of $\Q_1$ is induced by that of 
$\DF(N^i)$. Hence we can easily check that $\Q_1$, 
which can be identified with $\plin (\DF(\N))$, is isomorphic to 
$\DF(\cH(\N)) \cong \bigoplus_{i \in \ZZ} \DF(H^i(\N))[-i]$. 
We can prove the statement for $\ilin (\DG(\M))$ in the same way.  
\end{proof}

\section{Weakly Koszul Modules}

Let $B$ be a noetherian Koszul algebra (e.g. $B= A$ or $A^!$) 
with the graded Jacobson radical $\m$. For $M \in \grB$ and 
an integer $i$, $M_{\< i \>}$ denotes the submodule of $M$ generated by 
its degree $i$ component $M_i$.

\begin{prop}\label{my old result}
In the above situation, the following are equivalent. 
\begin{itemize}
\item[(1)] $M_{\< i \>}$ has a linear projective resolution for all $i$. 
\item[(2)] $H^i(\plin(M)) = 0$ for all $i \ne 0$. 
\item[(3)] All indecomposable summands of 
$\grm M$ have linear resolutions as $B \, (\cong \grm B)$-modules.  
\end{itemize}
\end{prop}

\begin{proof}
This result was proved in \cite[Proposition~4.9]{Y} under the assumption 
that $B$ is a polynomial ring. (R\"omer also proved this for a 
commutative  Koszul algebra. See \cite[Theorem~3.2.8]{R02}.)
In this proof, only the Koszul property of a polynomial ring 
is essential, and the proof also works in our case. But, to refer this, 
the reader should be careful with the following points. 

\smallskip

(a) In \cite{Y}, the grading of $\grm M$ is given by a different way. 
There, $(\grm M)_i = \m^i M /\m^{i+1} M$. It is easy to see that $\grm M$ 
has a linear resolution in this grading if and only if 
the condition (3) of the proposition is satisfied in our grading.  

\smallskip

(b) In the proof of \cite[Proposition~4.9]{Y},  
the regularity $\reg(N)$ of $N \in \grB$ is  an important tool. 
Unless $B$ is AS regular, one cannot define $\reg(N)$ 
using the local cohomologies of $N$. But if we set 
$\reg(N) := \sup \{ \, i+j \mid  \beta^{i,j} (N) \ne 0 \, \},$ 
then everything works well. It is not clear whether 
$\reg(N) < \infty$ for all $N \in \grB$ (c.f. \cite{Jor04}). 
But modules appearing in the argument similar to 
the proof of \cite[Proposition~4.9]{Y} have finite regularities. 

\smallskip

(c) In the proof of \cite[Proposition~4.9]{Y}, a few basic properties 
of the Castelnuovo-Mumford regularity (over a polynomial ring) are used. 
But $\reg(N)$ of $N \in \grB$ also has these properties, if 
we define $\reg(N)$ as (b). 
For example, if $N \in \grB$ satisfies $\dim_K N < \infty$, then 
$\reg(N) = \sigma(N)$. This can be proved by induction on 
$\dim_K N$. Using the short exact sequence $0 \to N_{\geq r} \to 
N \to N/N_{\geq r} \to 0$, we can also prove that $N_{\geq r}$ 
has an $r$ linear resolution if and only if $r \geq \reg(N)$ 
(see also Proposition~\ref{truncate}).  

\smallskip

(d) For the implication $(2) \Rightarrow (3)$, \cite{Y} refers another paper. 
But this implication can be proved by spectral sequence argument, 
since $\plin(M)$ can be constructed using a spectral sequence 
as we have seen in the previous section. 
\end{proof}

\begin{dfn}[{\cite{GM, MZ}}]
In the above situation, we say $M \in \grB$ is {\it weakly Koszul}, 
if it satisfies the equivalent conditions of Proposition~\ref{my old result}. 
\end{dfn}

\begin{rem}\label{wKos remarks}
(1) If $M \in \grB$ has a linear resolution, then it is weakly Koszul. 

(2) The notion of weakly Koszul modules was first introduced by 
Green and Martinez-Villa \cite{GM}. But they used the 
name ``strongly quasi Koszul modules". 
Weakly Koszul modules are also called ``componentwise 
linear modules" by some commutative algebraists (see \cite{HH}). 
\end{rem}

\begin{thm}\label{wKos} 
Let $0 \ne N \in \grAq$ and set $N' := \DAq(N)$.  
Then the following are equivalent. 
\begin{itemize}
\item[(1)] $N$ is weakly Koszul.  
\item[(2)] $H^i(\DFo (N'))$ has a $(-i)$-linear 
projective resolution for all $i$. 
\item[(3)] $\reg( \, \cH \circ \DFo (N') \,) = 0$. 
\item[(4)] $\reg( \, \cH \circ \DFo (N') \,) \leq 0$.
\end{itemize}
\end{thm}

\begin{proof}
Since $\iota(\cH \circ \DFo (N'))) \geq 0$ (i.e., 
$\iota(H^i( \DFo (N'))) \geq -i$ for all $i$), 
the equivalence among (2), (3) and (4) follows from 
Proposition~\ref{truncate}. So it suffices to prove 
$(1) \Leftrightarrow (4)$. Since ${\bf D}_{\Aqo}(\ilin(N')) \cong \plin(N)$, 
$N$ is weakly Koszul if and only if 
$H^i(\ilin( N'))=0$ for all $i > 0$. 
By Proposition~\ref{strand}, we have 
$$\ilin(N') = \ilin(\DGo \circ \DFo (N')) 
= \DGo \circ \cH \circ  \DFo (N').$$
Therefore, by Corollary~\ref{reg & DG},
$H^i(\ilin( N'))=0$ for all $i > 0$ 
if and only if the condition (4) holds.  
\end{proof}

\begin{rem}
Martinez-Villa and Zacharia proved that if $N$ is weakly Koszul then 
there is a filtration 
$$U_0 \subset U_1 \subset \cdots \subset U_p = N$$ 
such that $U_{i+1}/U_i$ has a linear resolution for each $i$ 
(see \cite[pp. 676--677]{MZ}). 
We can interpret this fact using Theorem~\ref{wKos} in our case. 

Let $N \in \grAq$ be a weakly Koszul module. Set $N' := \DAq(N)$ and 
$\L := \DFo(N')$.  Assume that $N$ does not have a linear 
resolution. Then $H^i(\L) \ne 0$ for several $i$. 
Set $n = \min \{ \, i \mid H^i(\L) \ne 0 \, \}$. 
Consider the truncation
$$\sigma_{> n} \L : \cdots \too 0 \too \operatorname{im}d^n 
\too T^{n+1} \too T^{n+2} \too \cdots$$ of $\L$. 
Then we have $H^i(\L) = H^i(\sigma_{> n} \L )$ for all $i > n$ and 
$H^i(\sigma_{> n} \L ) = 0$ for all $i \leq n$. 
We have a triangle 
\begin{equation}\label{triangle}
H^n(\L)[-n] \to \L \to \sigma_{> n} \L \to H^n(\L)[-n+1].
\end{equation} 
By Theorem~\ref{wKos}, $H^n(\L)[-n]$ has a 0-linear resolution.  
On the other hand, $$0 = \reg ( \, \cH(\sigma_{> n} \L) \, ) \geq  
\reg (\sigma_{> n} \L) \geq \iota(\sigma_{> n} \L) \geq 0.$$ 
Hence $\sigma_{> n} \L$ also has a 0-linear resolution.  Therefore, both 
${\bf D}_{\Aqo} \circ \DGo (\sigma_{> n} \L)$ and 
${\bf D}_{\Aqo} \circ \DGo ( H^n(\L)[-n]  )$ are acyclic complexes 
(that is, the $i^{\rm th}$ cohomology vanishes for all $i \ne 0$).  
Set $$U := H^0( \, {\bf D}_{\Aqo} \circ \DGo (\sigma_{> n} \L) \, ) 
\quad  \text{and} \quad 
V := H^0( \, {\bf D}_{\Aqo} \circ \DGo ( H^n(\L)[-n]  ) \, ).$$  
Since $N = {\bf D}_{\Aqo} \circ \DGo (\L)$, the triangle 
\eqref{triangle} induces a short exact sequence 
$0 \to U \to N \to V \to 0$ in $\grAq$. 
It is easy to see that $V$ has a linear resolution.  
Since $\cH \circ \DFo \circ \DAq (U) = 
\cH(\sigma_{> n} \L)$, $U$ is weakly Koszul by 
Theorem~\ref{wKos}. Repeating this procedure, we can get the expected 
filtration.  
\end{rem}

Let $N \in \grAq$ and $\cdots \stackrel{f_2}{\too} 
P^{-1} \stackrel{f_1}{\too} P^0 \stackrel{f_0}{\too} N \to 0$ 
its minimal projective resolution. For $i \geq 1$, 
we call $\Syz_i(N) := \ker (f_{i-1})$ 
the $i^{\rm th}$ {\it syzygy} of $N$. 
Note that $\Syz_i(N) = \image(f_i) = \coker(f_{i+1})$. 

By the original definition of a weakly Koszul module given in \cite{GM, MZ}, 
if $N \in \grAq$ is weakly Koszul then so is $\Omega_i(N)$ for all 
$i \geq 1$. 

\begin{dfn}[Herzog-R\"omer, \cite{R02}]
For $0 \ne N \in \grAq$, set
$$\lpd(N):= \inf\{ \, i \in \NN  \mid  \text{$\Syz_i(N)$ is 
weakly Koszul} \, \},$$
and call it the {\it linear part dominates} of $N$. 
\end{dfn}

Since $A \, (\cong (A^!)^!)$ 
is a noetherian ring of finite global dimension,   
$\lpd(N)$ is finite for all $N \in \grAq$ by \cite[Theorem~4.5]{MZ}.  

\begin{thm}\label{wKos=regH}
Let $N \in \grAq$ and set $N' := \DAq(N)$.   Then we have 
\begin{eqnarray*}
\lpd(N) &=& \reg( \, \cH \circ \DFo (N') \,)\\
&=& \max 
\{ \, \reg(H^i(\DFo(N') )) +i \mid i \in \ZZ \, \}.  
\end{eqnarray*}
\end{thm}

\begin{proof}
Note that $\P := {\bf D}_{\Aqo} \circ \DGo \circ \DFo(N')$ 
is a projective resolution of $N$, and 
$\Q := {\bf D}_{\Aqo} \circ \DGo (\DFo(N')_{\geq i})$ is the truncation 
$\cdots \to P^{-i-1} \to P^{-i} \to 0 \to \cdots $ 
of $\P$ for each $i \geq 1$. 
Hence we have $H^{j}(\Q) = 0$ for all $j \ne -i$ and 
there is a projective module $P$ such that 
$H^{-i}(\Q) \cong \Syz_i(N) \oplus P$. 
Since $P$ is weakly Koszul, $\Syz_i(N)$  is weakly 
Koszul if and only if so is $Q := H^{-i}(\Q)$. We have 
$$\plin (Q)[i] \cong \DAqo \circ \DGo \circ \cH 
(\DFo(N')_{\geq i}).$$
By Theorem~\ref{wKos}, $Q$ is weakly Koszul if and only if 
$\cH (\DFo(N')_{\geq i})$ has an $i$-linear resolution, that is, 
$H^j (\DFo(N')_{\geq i})$ has an 
$(i-j)$-linear resolution for all $j$. But there is 
some $L \in \grAqo$ such that $L = L_{i-j}$ and 
$H^j (\DFo(N')_{\geq i}) \cong H^j(\DFo(N'))_{\geq i-j} \oplus L$. 
Note that $L$ has an $(i-j)$-linear resolution. 
Therefore, $H^j (\DFo(N')_{\geq i})$ has an 
$(i-j)$-linear resolution if and only if so does 
$H^j(\DFo(N'))_{\geq i-j}$. Summing up the above facts, we have that 
$\Syz_i(N)$ is weakly Koszul if and only if 
$(\cH \circ \DFo(N'))_{\geq i}$ has an $i$-linear resolution.   
By Proposition~\ref{truncate}, the last condition is 
equivalent to the condition that $i \geq \reg (\cH \circ \DFo(N'))$. 
\end{proof}

\begin{rem}
Assume that $A$ is noetherian, Koszul, and has finite global dimension, 
but not necessarily AS regular. Then $A^!$ is a finite dimensional Koszul  
algebra, but not necessarily selfinjective. 
Even in this case, $\DG(\M)$ for $\M \in D^b(\grA)$ 
is a complex of injective $A^!$-modules, 
and the results in \S3 and Theorem~\ref{wKos=regH} also hold. 
But now we should set 
$\reg(\M) := \sup \{ \, i+j \mid  \beta^{i,j} (\M) \ne 0 \, \}$  
for $\M \in D^b(\grA)$ (local cohomology is not 
helpful to define the regularity). 
Since  $A$ is noetherian and has finite global 
dimension, we have $\reg(\M) < \infty$ for all $\M$.
In particular, we have $\lpd(N) < \infty$ for all $N \in \grAq$ 
(if $A$ is {\it right} noetherian) as proved in \cite[Theorem~4.5]{MZ}. 
\end{rem}

If $\lpd(N) \geq 1$ for some $N \in \grAq$, then 
$\sup \{ \, \lpd(T) \mid  T \in \grAq \, \} = \infty$. 
In fact, if $\Omega_{-i}(N)$ is the $i^{\rm th}$ {\it cosyzygy} of $N$ 
(since $A^!$ is selfinjective, we can consider cosyzygies), 
then $\lpd(\Omega_{-i}(N)) > i$. 
But Herzog and R\"omer proved that if $J$ is a monomial ideal 
of an exterior algebra $E =  \bigwedge \< y_1, \ldots, y_d \>$ 
then $\lpd (E/J) \leq d-1$ (c.f. \cite[\S3.3]{R02}).  
We will refine their results using Theorem~\ref{wKos=regH}. 

\medskip

In the sequel, we regard the polynomial ring 
$S=K[x_1, \ldots, x_d]$, $d \geq 1$, as an $\NN^d$-graded ring 
with $\deg x_i = (0, \ldots, 0,1,0, \ldots, 0)$ where 1 is at the $i^{\rm th}$ 
position. Similarly, the exterior algebra $E = S^!= \bigwedge
\<y_1, \ldots, y_d \>$ is also an  $\NN^d$-graded ring. 
Let $\MMSn$ be the category of $\ZZ^d$-graded $S$-modules and their 
degree preserving $S$-homomorphisms, 
and $\mmSn$ its full subcategory consisting 
of finitely generated modules. We have a similar category $\mmEn$ for $E$. 
For $S$-modules and graded $E$-modules, we do not have to 
distinguish left modules from right modules. Since $\ZZ^d$-graded modules 
can be regarded as $\ZZ$-graded modules in the natural way, 
we can discuss $\reg(\M)$ for $\M \in D^b(\mmSn)$ and 
$\lpd(N)$ for $N \in \mmEn$. 

Note that $\DE(-) = \bigoplus_{\ba \in \ZZ^d} 
\Hom_{\mmEn}(-,E(\ba))$ gives an exact duality functor 
from $\mmEn$ to itself. Sometimes, we simply denote $\DE(N)$ by $N'$. 
Set $\b1 := (1, 1, \ldots, 1) \in \ZZ^d$. Then 
$\D_S := S(-\b1)[d] \in D^b(\mmSn)$ is a $\ZZ^d$-graded 
normalized dualizing complex and $\DS(-):= \bR \uHom_S(-,\D_S) 
= \bigoplus_{\ba \in \ZZ^d} \bR \Hom_{\MMSn} (-, \D_S(\ba))$ 
gives a duality functor from $D^b(\mmSn)$ to itself. 
As shown in \cite[Theorem~4.1]{Y5}, we have the $\ZZ^d$-graded 
Koszul duality functors $\DFn$ and $\DGn$ giving an equivalence 
$D^b(\mmSn) \cong D^b(\mmEn)$. These functors are defined in the same 
way as in the $\ZZ$-graded case. 

For $\ba = (a_1, \ldots, a_d)\in \ZZ^d$, 
set $\supp (\ba) := \{i  \mid a_i > 0\} \subset [d] := \{1, \ldots, d \}$. 
We say  $\ba \in \ZZ^d$ is {\it squarefree} if $a_i= 0,1$ for all $i \in [d]$. 
When $\ba \in \ZZ^d$ is squarefree, we sometimes identify $\ba$ with  
$\supp (\ba)$. For example, if $F \subset [d]$, then 
$S(-F)$ means the free module $S(-\ba)$, where $\ba \in \NN^d$ 
is the squarefree vector with $\supp(\ba) = F$. 

\begin{dfn}[\cite{Y}]
We say $M \in \mmSn$ is {\it squarefree}, 
if $M$ has a presentation of the form 
$$\bigoplus_{F \subset [d]} S(-F)^{m_F} \to 
\bigoplus_{F \subset [d]} S(-F)^{n_F} \to M \to 0$$
for some $m_F, n_F \in \NN$. 
\end{dfn}

The above definition seems different from 
the original one given in \cite{Y}, but they coincide. 
Stanley-Reisner rings (that is, the quotient rings of $S$ by 
squarefree monomial ideals) and many modules related to them 
are squarefree. 
Here we summarize basic properties of squarefree 
modules. See \cite{Y, Y5} for further information. 
Let $\Sq$ be the full subcategory of $\mmSn$ consisting of 
squarefree modules. 
Then $\Sq$ is closed under kernels, cokernels, and 
extensions in $\mmSn$. Thus $\Sq$ is an abelian category. Moreover, we have  
$D^b(\Sq) \cong D^b_{\Sq}(\MMSn)$. 
If $M$ is squarefree, then each term in a 
$\ZZ^d$-graded minimal free resolution of $M$ is of the form 
$\bigoplus_{F \subset [d]} S(-F)^{n_F}$. Hence we have $\reg(M) \leq  d$. 
Moreover, $\reg(M) = d$ if and only if $M$ has a summand 
which is isomorphic to $S(-\b1)$. 

\begin{dfn}[R\"omer~\cite{R0}]\label{sqE}
We say $N \in \mmEn$ is 
{\it squarefree}, if $N = \bigoplus_{F \subset [d]} N_F$ 
(i.e., if $\ba \in \ZZ^d$ is not squarefree, then $N_\ba = 0$). 
\end{dfn}

A monomial ideal of $E$ is always a squarefree $E$-module. 
Let $\SqE$ be the full subcategory of $\mmEn$ consisting of 
squarefree modules. Then $\SqE$ is an abelian category with 
$D^b(\SqE) \cong D^b_{\SqE}(\mmEn)$.  
If $N$ is a squarefree $E$-module, 
then so is $\DE(N)$. That is, $\DE$ gives an exact duality functor 
from $\SqE$ to itself.  We have functors 
$\cS: \SqE \to \Sq$ and $\cE: \Sq \to \SqE$ giving an 
equivalence $\Sq \cong \SqE$. Here 
$\cS(N)_F = N_F$ for $N \in \SqE$ and $F \subset [d]$, 
and the multiplication map $\cS(N)_F \ni z \mapsto x_i z \in 
\cS(N)_{F \cup \{ i \}}$ for $i \not \in F$ is given by 
$\cS(N)_F =N_F \ni z \mapsto (-1)^{\alpha(i, F)} y_i z 
\in N_{F \cup \{ i \}}= \cS(N)_{F \cup \{ i \}}$, where 
$\alpha(i,F) =  \# \{ \, j \in F \mid j < i \, \}$. 
See \cite{R0, Y5} for detail. 
Since a free module $E(\ba)$ is {\it not} squarefree 
unless $\ba = 0$, the syzygies of a squarefree $E$-module are 
{\it not} squarefree. 

\begin{prop}[{Herzog-R\"omer, \cite[Corollary~3.3.5]{R02}}]
\label{Tim} If $N$ is a squarefree $E$-module 
(e.g., $N = E/J$ for a monomial ideal $J$), 
then we have $\lpd(N) \leq d-1$. 
\end{prop}

This result easily follows from  
Theorem~\ref{wKos=regH} and the fact that 
$H^i(\DFn(N'))(-\b1)$ is a squarefree $S$-module for all $i$  
and $H^i(\DFn(N')) = 0$ unless $0 \leq i \leq d$. 
(Recall the remark on the regularity of squarefree modules 
given before Definition~\ref{sqE}, and note that 
$M:=H^d(\DFn(N'))(-\b1)$ is generated by $M_0$). 

We also remark that \cite[Corollary~3.3.5]{R02} just states 
that $\lpd (N) \leq d$. 
But their argument actually proves that $\lpd (N) \leq d-1$.   
In fact, they showed that $$\lpd(N) \leq \pdim_S \cS(N).$$ 
But, if $\pdim_S \cS(N) = d$ then $\cS(N)$ has a summand   
which is isomorphic to $K = S/(x_1, \ldots, x_d)$ and hence 
$N$ has a summand which is isomorphic to 
$K = E/(y_1, \ldots, y_d)$. But $K \in \SqE$ has 
a linear resolution and irrelevant to $\lpd(N)$.  
\medskip

To refine Proposition~\ref{Tim}, we need further properties of 
squarefree modules.

If $\M \in D^b(\Sq)$, then 
$\uExt^i_S(\M, \D_S)$ is squarefree for all $i$. Hence 
$\D_S$ gives a duality functor on $D^b(\Sq)$. On the other hand, 
$\bA := \cS \circ \DE \circ \cE$ is an exact 
duality functor on $\Sq$ and it induces a duality functor on $D^b(\Sq)$.  
Miller \cite[Corollary~4.21]{Mil} and 
R\"omer \cite[Corollary~3.7]{R1} proved that 
$\reg(\bA(M)) = \pdim_S M$ for all $M \in \Sq$. 
I generalized this equation to a complex $\M \in D^b(\Sq)$ in 
\cite[Corollary~2.10]{Y7}.
 
\begin{lem}\label{sqf formula}
Let $N \in \SqE$ and set $N' := \DE(N)$. Then we have 
\begin{equation}\label{1st}
\reg(H^i(\DFn(N'))) = - \depth_S 
( \, \Ext_S^{d-i}( \CS(N'), \, S ) \, )  
\end{equation} and
\begin{equation}\label{2nd}
\lpd(N) = \max \{ \, i - \depth_S 
( \, \Ext_S^{d-i}( \CS(N'), \, S ) \, )  \mid  
0 \leq i \leq d  \, \}.
\end{equation}
Here we set the depth of the 0 module to be $+\infty$.
\end{lem}

If $M :=  \Ext_S^{d-i}( \CS(N'), \, S ) \ne 0$,  
then $\depth_S M \leq \dim_S M \leq i$. 
Therefore all members in the set of the right side of 
\eqref{2nd} are non-negative or $- \infty$.

\begin{proof}
By Theorem~\ref{wKos=regH}, \eqref{2nd} follows from \eqref{1st}. 
So it suffices to show \eqref{1st}. 
By \cite[Proposition~4.3]{Y5}, we have $\DFn(N') \cong 
(\bA \circ \DS \circ \cS (N'))(\b1)$. (The degree shifting 
``$(\b1)$" does not occur in \cite[Proposition~4.3]{Y5}. 
But $E$ is a negatively graded ring there, 
and we need the degree shifting in the present convention.) 
Since $\bA$ is exact, we have 
\begin{eqnarray*}
H^i(\DFn(N'))) \cong H^i( \, \bA \circ \DS \circ \cS (N') \, ) (\b1) 
&\cong& \bA ( \, H^{-i}(\DS \circ \cS (N')) \, ) (\b1)\\ 
&=& \bA ( \, \uExt_S^{-i} (\cS(N'), \, \D_S) \, ) (\b1).
\end{eqnarray*}
Recall that $\reg(\bA(M)) = \pdim_S M$ for $M \in \Sq$. 
On the other hand, since $M$ is finitely generated, 
the underlying module of $\uExt_S^{-i}(M,\D_S)$ is 
isomorphic to $\Ext_S^{d-i}(M, S)$. 
So \eqref{1st} follows from these facts and 
the Auslander-Buchsbaum formula. 
\end{proof}

\begin{cor}\label{comp lin}
For $N \in \SqE$, $N$ is weakly Koszul (over $E$) 
if and only if $\cS(N)$ is weakly Koszul (over $S$). 
\end{cor}

In \cite[Corollary~1.3]{R1}, it was proved that $N$ has a linear 
resolution if and only if so does $\cS(N)$. 
Corollary~\ref{comp lin} also follows from this fact and 
(the squarefree module version of) \cite[Proposition~1.5]{HH}. 

\begin{proof}
We say $M \in \grS$ is {\it sequentially Cohen-Macaulay}, 
if for each $i$ $\Ext_S^i(M,S)$ is either the zero module or 
a Cohen-Macaulay module of dimension $d-i$ 
(c.f. \cite[III. Theorem~2.11]{St}). 
By Lemma~\ref{sqf formula}, $N$ is weakly Koszul 
if and only if $\cS(N')$ ($\cong \bA \circ \cS (N)$) 
is sequentially Cohen-Macaulay.  
By \cite[Theorem~4.5]{R1}, the latter condition holds  
if and only if $\cS(N)$ is weakly Koszul. 
\end{proof}

Many examples of squarefree monomial ideals of $S$ 
which are weakly Koszul (dually, Stanley-Reisner rings which are 
sequentially Cohen-Macaulay) are known. 
So we can obtain many weakly Koszul monomial 
ideals of $E$ using Corollary~\ref{comp lin}. 

\begin{prop}
For an integer $i$ with $1 \leq i \leq d-1$,  
there is a squarefree $E$-module $N$ such that 
$\lpd N = \pdim_S \cS(N) = i$. 
In particular, the inequality of 
Proposition~\ref{Tim} is optimal. 
\end{prop}

\begin{proof}
Let $M$ be the $\ZZ^d$-graded $i^{\rm th}$  syzygy of 
$K = S/\m$. Note that $M$ is squarefree. We can easily 
check that $N:=\DE \circ \cE(M) \in \SqE$ satisfies the expected 
condition.   In fact, $\pdim_S \cS(N) 
= \pdim_S \bA (M) = \reg M = i$.  On the other hand, 
since  $\Ext^{d-i}_S(\cS(N'), S) = \Ext^{d-i}_S(M, S) = K$,  
$\Ext^j_S(\cS(N'), S)=0$ for all $j \ne d-i, 0$, and 
$\depth_S (\Hom_S(\cS(N'), S)) = d-i+1$, we have $\lpd N = i$. 
\end{proof}

The above result also says that the inequality 
$\lpd(N) \leq \pdim_S \cS(N)$ of \cite[Corollary~3.3.5]{R02}
is also optimal. 
But for a monomial ideal $J \subset E$, the situation is different.  

\begin{prop}\label{d-2}
If $d \geq 3$, then we have 
$\lpd (E/J) \leq d-2$ for a monomial ideal $J$ of $E$. 
\end{prop}

\begin{proof}
If $d=3$, then easy computation shows that 
any squarefree monomial ideal $I \subset S$ is weakly Koszul. 
Hence $J$ is weakly Koszul by Corollary~\ref{comp lin}. 
So we may assume that $d \geq 4$. 

Note that $\bA \circ \cS(E/J)$ is isomorphic to a squarefree monomial 
ideal of $S$. We denote it by $I$. 
By Lemma~\ref{sqf formula}, 
it suffices to show that $\depth_S (\Hom_S(I,S)) \geq 2$ 
and $\depth_S (\Ext_S^1(I,S)) \geq 1$. 
Recall that $\Hom_S(I,S)$ satisfies  
Serre's condition $(S_2)$, hence its depth is 
at least 2. Since $\Ext_S^1(I,S) \cong 
\Ext_S^2(S/I,S)$, it suffices to prove 
that $\depth_S (\Ext_S^2(S/I, S)) \geq 1$. 

If $\cod (I) > 2$, then we have $\Ext_S^2(S/I,S) = 0$.   
If $\cod(I) = 2$, then $\Ext_S^2(S/I,S)$ satisfies $(S_2)$ as an 
$S/I$-module and $\depth_S \Ext_S^2(S/I,S) \geq \min \{ \, 2, 
\dim (S/I) \, \} \geq 2$. So we may assume that $\cod(I)=1$. 
If the heights of all associated primes of $I$ are 1, 
then $I$ is a principal ideal and $\Ext^i_S(S/I, S) = 0$ 
for all $i \ne 1$. So we may assume that $I$ has an prime of larger 
height. 
Then we have ideals $I_1$ and $I_2$ of $S$ such that $I = I_1 \cap I_2$ 
and the heights of any associated prime of $I_1$ 
(resp. $I_2$) is 1 (at least 2).
Since $I$ is a radical ideal, we have $\cod (I_1 + I_2) \geq 3$. 
Hence $\Ext^2_S(S/(I_1 + I_2), S) = 0$ and 
$\Ext^3_S(S/(I_1 + I_2), S)$ is either the zero module or 
satisfies $(S_2)$ as an $S/(I_1+I_2)$-module.  
In particular, if $\Ext_S^3(S/(I_1 + I_2), S) \ne 0$ 
(equivalently, if $\dim (S/(I_1 + I_2)) = d-3$) then  
$\depth_S (\Ext^3_S(S/(I_1 + I_2), S)) 
\geq \min\{ 2, d-3 \} \geq 1$. Note that  
$\depth_S ( \Ext_S^2(S/I_2, S) ) \geq 2$. 
From the short exact sequence 
$$0 \to S/I \to S/I_1 \oplus S/I_2 \to S/(I_1 + I_2) 
\to 0$$ and the above argument, we have the exact sequence  
\begin{equation}\label{last}
0 \to \Ext^2_S(S/I_2, S) \to \Ext^2_S(S/I, S) \to 
\Ext^3_S(S/(I_1 + I_2), S).
\end{equation}  
We have $\depth_S (\Ext^2_S(S/I,S)) \geq 1$ by \eqref{last}, 
since the modules beside this module have positive depth. 
\end{proof}

\section*{Acknowledgments}
This research  started during my stay at Instituto de Matem\'aticas 
de la UNAM (Morelia, Mexico) in the summer of 2004. 
I am grateful to JSPS and CONACYT for financial support, 
and the institute for warm hospitality.  
Special thank are due to Professor Robert Martinez-Villa for 
stimulating discussion and encouragement. 
I also thanks Professor Izuru Mori for valuable information on AS regular 
algebras.


\begin{thebibliography}{99}
\bibitem{ATV} M. Artin,  J. Tate,  and M. Van den Bergh,  
Some algebras associated to automorphisms of elliptic curves, in: 
The Grothendieck Festschrift, Vol. I, 33--85, 
Progr. Math., 86, Birkhauser Boston, Boston, MA, 1990. 

\bibitem{BGS}
A. Beilinson,  V. Ginzburg and   W. Soergel,  
Koszul duality patterns in representation theory, 
J. Amer. Math. Soc. 9 (1996), 473--527. 


\bibitem{EFS} D. Eisenbud, G. Fl\o ystad  and F.-O. Schreyer, 
Sheaf cohomology and free resolutions over exterior algebra, 
Trans. Amer. Math. Soc. 355 (2003), 4397-4426.  


\bibitem{EG} D. Eisenbud and S. Goto, 
Linear free resolutions and minimal multiplicity, 
{\it J. Algebra} 88 (1984), 89--133.

\bibitem{GM}
E.L. Green, R. Martinez-Villa, 
Koszul and Yoneda algebras, in: Representation Theory of Algebras 
(Cocoyoc, 1994), CMS Conf. Proc., Vol. 18, 
Amer. Math. Soc., 1996, pp. 247--297. 

\bibitem{GMT} 
J.Y. Guo, R. Martinez-Villa, and M. Takane, 
Koszul generalized Auslander regular algebras, in: 
Algebras and modules, II (Geiranger, 1996),  
CMS Conf. Proc., 24, Amer. Math. Soc., 1998, pp. 263--283. 

\bibitem{HH} J. Herzog and T. Hibi, 
Componentwise linear ideals, Nagoya Math. J. 153 
(1999), 141--153. 


\bibitem{HI} J. Herzog and S. Iyengar, Koszul modules, to appear in 
J. Pure Appl. Algebra.  


\bibitem{Jor99}P. J\o rgensen, 
Non-commutative Castelnuovo-Mumford regularity, 
Math. Proc. Cambridge Philos. Soc. 125 (1999) 203-221.

\bibitem{Jor04} P. J\o rgensen, 
Linear free resolutions over non-commutative algebras, 
Compos. Math. 140 (2004), 1053--1058.


\bibitem{M98} R. Martinez-Villa, 
Serre duality for generalized Auslander regular algebras, In: 
Trends in the representation theory of finite-dimensional 
algebras (Seattle, WA, 1997),  Contemp. Math., 229,
Amer. Math. Soc., 1998, pp. 237--263.  

\bibitem{M99}
R. Martinez-Villa,  
Graded, selfinjective, and Koszul algebras. 
J. Algebra 215 (1999),  34--72.

\bibitem{MZ} R. Martinez-Villa and  D. Zacharia,  
Approximations with modules having linear resolutions, 
J. Algebra 266 (2003), 671--697.

\bibitem{Mil} E. Miller, 
The Alexander duality functors and local duality with monomial support, 
J. Algebra. {\bf 231} (2000), 180--234. 

\bibitem{Mori} I. Mori,  
Rationality of the Poincare series for Koszul algebras, 
J. Algebra 276 (2004), 602--624.

\bibitem{R0} T. R\"omer,  Generalized Alexander duality and applications, 
Osaka J. Math. 38 (2001), 469--485.  

\bibitem{R1} T. R\"omer,  Cohen-Macaulayness and squarefree modules, 
Manuscripta Math., 104 (2001), 39--48.

\bibitem{R02} T. R\"omer, 
 On minimal graded free resolutions, Thesis, University of Essen, 2001. 
 
\bibitem{St} R. Stanley,  
Combinatorics and commutative algebra, 2nd ed. Birkh\"auser 1996.   


\bibitem{Y} K. Yanagawa, 
Alexander duality for Stanley-Reisner rings and squarefree 
$\NN^n$-graded modules, J. Algebra 225 (2000), 630--645.  

\bibitem{Y5} K. Yanagawa, 
Derived category of squarefree modules and local cohomology with monomial 
ideal support, J. Math. Soc. Japan 56 (2004) 289--308. 

\bibitem{Y7} K. Yanagawa, 
BGG correspondence and R\"omer's theorem of an exterior algebra, 
to appear in Algebras and Representation Theory. (math.AC/0402406) 


\bibitem{Ye} A. Yekutieli, 
Dualizing complexes over noncommutative graded algebras, 
J. Algebra 153 (1992), 41--84.  
\end{thebibliography}
\end{document}